\documentclass[12pt]{article}
\usepackage{amssymb}
\begin{document}

\def\ad{\mbox{ad}\,}
\def\ker{\mbox{Ker}\,}
\def\vev#1{\langle #1 \rangle}
\newcommand{\sect}[1]{\setcounter{equation}{0}\section{#1}}
\renewcommand{\theequation}{\thesection.\arabic{equation}}
\newcommand{\subsect}[1]{\setcounter{equation}{0}\subsection{#1}}
\renewcommand{\theequation}{\thesubsection.\arabic{equation}}
\newcommand{\be}{\begin{equation}}
\newcommand{\ee}{\end{equation}}
\newcommand{\bea}{\begin{eqnarray}}
\newcommand{\eea}{\end{eqnarray}}

\title{Noncommutative Geometry of the $h$-deformed Quantum Plane}
\author{S. Cho\\ Department of Physics, Semyung University, \\
        Chechon, Chungbuk 390 - 230, Korea
    \and
        J. Madore\\
        Laboratoire de Physique Th\'eorique et Hautes
        Energies\thanks{Laboratoire associ\'e au CNRS, URA D0063}\\
        Universit\'e de Paris-Sud, B\^at. 211, F-91405 Orsay, France
    \and
        K.S. Park \\
        Department of Mathematics, Keimyung University, \\
        Taegu 705 - 701, Korea
        }


\maketitle

\begin{abstract}
The $h$-deformed quantum plane is a counterpart of the $q$-deformed one
in the set of quantum planes which are covariant under those quantum
deformations of $GL(2)$ which admit a central determinant.  We have
investigated the noncommutative geometry of the $h$-deformed quantum
plane. There is a 2-parameter family of torsion-free linear connections,
a 1-parameter sub-family of which are compatible with a skew-symmetric
non-degenerate bilinear map.  The skew-symmetric map resembles a
symplectic 2-form and induces a metric.  It is also shown that the
extended $h$-deformed quantum plane is a noncommutative version of the
Poincar\'e half-plane, a surface of constant negative Gaussian
curvature.
\end{abstract}



\vfill
\noindent
Semyung Preprint SMHEP 97/8
\newpage

\parskip 4pt plus2pt minus2pt

\sect{ Introduction }

Quantum planes are simple examples of quantum spaces and have been
studied intensively by many authors in the past years. They arise as
deformations of planes on which quantum groups act covariantly. For
references to the literature we refer to the recent monographs by Chari
\& Pressley~\cite{ChaPre94} and by Majid~\cite{Maj95}.
One of the quantum planes, referred to as the $q$-deformed quantum plane
or the Manin plane~\cite{Man88}, is defined as the associative algebra
generated by two noncommuting elements (`coordinates') $x$ and $y$ such
that
$$
xy = q yx.
$$
The quantum group $GL_{q}(2) $ is the symmetry group of the $q$-quantum
plane. Another quantum plane, called the $h$-deformed quantum
plane~\cite{DemManMukZhd90, Man91}, is defined as the associative
algebra generated by two noncommuting elements $x$ and $y$ such that
$$
xy - yx = hy^2.
$$
The quantum group $GL_{h}(2)$ is the symmetry group of the $h$-quantum
plane.  These two quantum planes are the only deformations of the
ordinary plane which are covariant under the quantum deformations of
$GL(2)$ which admit a central determinant since up to isomorphism
$GL_{q}(2)$ and $GL_{h}(2)$ are the only two such deformed quantum
groups~\cite{Kup92}.  The $h$-deformation can be seen as a singular
contraction of a $q$-deformation~\cite{AghKhoSha95}. More precisely, a
class of similarity transformations of the $R$-matrices associated to
$q$-deformations can be introduced such that the $q \rightarrow 1$ limit
gives explicit $R$-matrices for the $h$-deformations~\cite{AbdChaCha97}.
Although the transformation matrix is itself singular in the limit, the
construction is well-defined.

As usual in noncommutative geometry~\cite{Con94,Mad95} quantum planes
have over them many differential calculi $\Omega^*({\cal A})$. The
commutation relations in $\Omega^1({\cal A})$ must be consistent with
the commutation relations of the algebra but this condition is not
enough to uniquely define the calculus. There is however a particularly
interesting calculus known as the Wess-Zumino calculus~\cite{PusWor89,
WesZum90} which is covariant under the co-action of the $q$-deformed
quantum groups. There is similarly a calculus over the $h$-deformed
quantum plane which is covariant under the co-action of the $h$-deformed
quantum groups~\cite{Agh93}.  Moreover, general definitions have been
proposed recently of a linear connection and a metric within the context
of noncommutative geometry in general and for quantum spaces in
particular. Using these tools, we shall here investigate the Riemannian
geometry of the $h$-deformed quantum plane.  It turns out that the
$h$-deformed quantum plane has more interesting geometrical properties
than the $q$-deformed one.

In Section~2 we give a review of the definition of what we call the
`Stehbein' formalism~\cite{DimMad96,MadMou96} and of a definition of a
linear connection~\cite{DubMadMasMou95,Mou95,Mad95,DubMadMasMou96} which
has been used in noncommutative geometry.  In Section~3, a 2-parameter
family of torsion-free linear connections is constructed on the
$h$-deformed quantum plane. The existence of a 2-parameter family of
torsion-free linear connections is shown to be quite general even within
the set of 2-parameter $h$-deformed quantum planes with an
appropriate supplementary condition between deforming parameters.
Moreover, there is a skew-symmetric non-degenerate bilinear map with
which a 1-parameter sub-family of linear connections are compatible.  We
shall also show that the skew-symmetric map resembles the symplectic
2-form of an ordinary manifold and induces a metric and the skew
derivatives~\cite{WesZum90} of the $h$-deformed quantum plane.  We shall
compare the results of the $h$-deformed quantum plane with those of the
$q$-deformed one.  In Section~4, we shall investigate the geometry of
the extended $h$-deformed quantum plane. It turns out that the extended
$h$-deformed quantum plane has a unique metric-compatible torsion-free
linear connection; it is a noncommutative version of the Poincar\'e
half-plane, a surface of constant negative Gaussian curvature. This can
be shown explicitly by a change of generators.

\sect{ Metric-compatible Linear connections }

\subsect{ Linear connections }

Let ${\cal A}$ be an associative algebra with the identity 1.  Let
$(\Omega^*_{u}, d_{u})$ be the universal differential calculus over
${\cal A}$. Then every other differential calculus over ${\cal A}$ can
be obtained as a quotient of it.  We suppose that there exists a
bimodule of 1-forms $\Omega^1$ and a map $d$ of ${\cal A}$ into
$\Omega^1$.  Then we can find an ${\cal A}$-bimodule homomorphism
$\phi_{1} : \Omega^1_{u} \rightarrow \Omega^1 $ such that
$\phi_{1} \circ d_{u} = d $.  For integers $n \geq 2$, $\Omega^n $ is
defined to be the quotient space
\be
\Omega^n \equiv \frac{\Omega^n_{u}}{\vev{d_{u}(\ker \phi_{n-1})}},
\ee
where $\phi_{n-1}$ is the projection map from $\Omega^{n-1}_{u}$ to
$\Omega^{n-1}$ and $\vev{d_{u}(\ker \phi_{n-1})}$ is the
${\cal A}$-bimodule generated by $d_{u}(\ker \phi_{n-1})$.  This
construction can be summarized in the following commutative diagram
\be
\def\normalbaselines{\baselineskip=18pt}
\matrix{
{\cal A} &\buildrel d_u \over \longrightarrow &\Omega_u^1
         &\buildrel d_u \over \longrightarrow &\Omega_u^2
         &\buildrel d_u \over \longrightarrow &\cdots                \cr
\parallel&&\phi_1 \downarrow \phantom{\phi_1}
         &&\phi_2 \downarrow \phantom{\phi_2}                        \cr
{\cal A} &\buildrel d \over \longrightarrow &\Omega^1
         &\buildrel d \over \longrightarrow &\Omega^2
         &\buildrel d \over \longrightarrow &\cdots .                \cr
}
\def\normalbaselines{\baselineskip=12pt}
\ee
We shall be mainly interested here in the bimodules $\Omega^1$ and
$\Omega^2$.  Since
\be
\Omega^2_u  = \Omega^1_u  \otimes_{\cal A} \Omega^1_u
\ee
there is  an exact sequence of ${\cal A}$-bimodules
\be
0 \rightarrow {\cal K} \hookrightarrow
\Omega^1 \otimes_{\cal A}\Omega^1
\buildrel \pi \over \rightarrow \Omega^2 \rightarrow 0            \label{eq:6}
\ee
where
\be
{\cal K} =
(\phi_1 \otimes \phi_1) (\vev{d_u \ker \phi_1}) =
(\phi_1 \otimes \phi_1) ({d_u \ker \phi_1}).
\ee

The definition of linear connection we
use~\cite{DubMadMasMou95,Mou95,Mad95,DubMadMasMou96} makes full use of
the bimodule structure of $\Omega^1$. It is defined to be a map
\be
D: \Omega^1 \longrightarrow \Omega^1 \otimes_{{\cal A}} \Omega^1  \label{eq:8}
\ee
satisfying the two Leibniz rules
\bea
D (f\xi)  & = &df \otimes \xi + f D \xi, \\                      \label{eq:9}
D (\xi f) & = &\sigma(\xi \otimes df) + (D \xi)f,               \label{eq:10}
\eea
where $f \in {\cal A}$, $\xi \in \Omega^1$ and $\sigma$ is a map from
$\Omega^1 \otimes_{A} \Omega^{1}$ to $\Omega^1 \otimes_{A} \Omega^1$
which generalizes the permutation.  For the consistency of the
definition of a linear connection, the $\sigma$ must be assumed to be
${\cal A}$-bilinear. That is, we must have, for $f \in A$ and
$\xi, \eta \in \Omega^1$,
\be
\sigma (f \xi \otimes \eta ) =  f \sigma ( \xi \otimes \eta ),
\hspace{1cm}
\sigma ( \xi \otimes \eta f ) =  \sigma ( \xi \otimes \eta ) f.  \label{eq:11}
\ee

The map $\Theta : \Omega^1 \rightarrow \Omega^2 $ defined by
$\Theta = d - \pi \circ D$ is the torsion of the linear connection $D$.
It is ${\cal A}$-bilinear only if $\sigma$ is assumed to satisfy the
condition~\cite{DubMadMasMou95}.
\be
\pi \circ (\sigma + 1) = 0.                                      \label{eq:12}
\ee

A linear connection $D$ can be extended to a linear map
\be
D :\Omega^1 \otimes_{A} \Omega^1 \longrightarrow
\Omega^1 \otimes_{A} \Omega^1 \otimes_{A} \Omega^1
\ee
satisfying
\be
D(\xi \otimes \eta )  =  D \xi \otimes \eta
+ \sigma_{12}(\xi \otimes D \eta)                                \label{eq:13}
\ee
for $\xi, \eta \in \Omega^1 $, where $\sigma_{12} = \sigma \otimes 1$.

An ${\cal A}$-bilinear map
\be
\Omega^1 \otimes _{{\cal A}} \Omega^1 \buildrel g \over
\rightarrow {\cal A}
\ee
is called non-degenerate whenever $g(\xi \otimes \eta) = 0 $ for all
$\eta \in \Omega^1$ implies that $\xi = 0$ and
$g(\xi \otimes \eta) = 0$ for all $\xi \in \Omega^1$ implies that
$\eta = 0$. A metric is a non-degenerate ${\cal A}$-bilinear map.
A metric is called symmetric (skew-symmetric)
if $g \circ \sigma = g$ ($g \circ \sigma = - g$).

A linear connection $D$ is said to be compatible with a metric
$g$ if the condition
\be
(1\otimes g) \circ D = d \circ g.                                \label{eq:14}
\ee
is satisfied

The  curvature is defined to be the map
\be
\pi_{12} D^2 : \Omega^1 \rightarrow \Omega^2 \otimes _{{\cal A}} \Omega^1,
\ee
where $\pi_{12} = \pi \otimes 1$.  This map is left ${\cal A}$-linear
but it is not in general right ${\cal A}$-linear~\cite{DubMadMasMou96}.
There is at the moment no general consensus of the correct definition of
the curvature respecting the bimodule structure of the linear connection
but since we are primarily interested in the first-order effects in the
commutative limit, we can identify the curvature with the operator
$\pi_{12}D^2$.

We define the Ricci map
\be
\Omega^1 \buildrel  \mbox{Ric} \over \longrightarrow   \Omega^1  \label{eq:15}
\ee
by $\mbox{Ric} = - (1\otimes g) D^2$.

\subsect{ The Stehbein formalism }

To initiate the construction in the previous subsection, we suppose that
the algebra ${\cal A}$ is noncommutative and define the bimodule of
1-forms using a set of inner derivations~\cite{DimMad96}. For each
positive integers $n$ let $\lambda_i$ be a set of $n$ linearly
independent elements of ${\cal A}$ and define the derivations by
\be
e_i = \ad \lambda_i.                                             \label{eq:16}
\ee
For any $f \in {\cal A}$, we define the 1-form $df$ by
\be
df(e_i) = e_i f = [\lambda_i, f].                                \label{eq:17}
\ee
The $\Omega^{1}$ is then defined to be the ${\cal A}$-bimodule
$\vev{d{\cal A}}$ generated by the image of $d$. Any element of
$\Omega^1$ is the sum of elements of the form $fdg$ or, equivalently
using the Leibniz rule, of the form $(df)g$.  We define
\be
(fdg) (e_i) = fe_ig, \hspace{2cm} ((dg)f)(e_i) = (e_i g)f.       \label{eq:18}
\ee
We suppose that there exists a set of $n$ elements $\theta^i$ of
$\Omega^1$, called~\cite{DimMad96,MadMou96} a `frame' or `Stehbein' as the
noncommutative equivalent of a `moving frame' or $n$-bein, such that
\be
\theta^i (e_j ) = \delta^i_j.                                    \label{eq:19}
\ee
Then it follows easily that $\theta^i$ commute with the elements
$f \in {\cal A }$,
\be
f\theta^i = \theta^i f.                                          \label{eq:20}
\ee
and that $\Omega^1$ is free of rank $n$ as a left or right module.
Hence the exact sequence in Equation~ (\ref{eq:6}) splits.  Let $\jmath$ be
the splitting map and write~\cite{MadMou96}
\be
\jmath \circ \pi (\theta^i \otimes \theta^j ) =
 P^{ij}{}_{kl}\,\, \theta^k \otimes \theta^l.                    \label{eq:23}
\ee
The coefficients $P^{ij}{}_{kl}$ depend on the map $\jmath $ and belong
to the center ${\cal Z(A)}$ of ${\cal A}$. Since $\pi $ is a projection
we have
\be
P^{ij}{}_{mn}P^{mn}{}_{kl} = P^{ij}{}_{kl}
\ee
and the product $\theta^i \theta^j $ satisfies the condition
\be
\theta^i \theta^j = P^{ij}{}_{kl} \theta^k \theta^l.             \label{eq:25}
\ee
If we define $\theta = - \lambda_i\theta^i $, then it follows that
\be
df = -[\theta, f]                                                \label{eq:21}
\ee
and thus, as an ${\cal A}$-bimodule, the one element
$\theta $ generates $\Omega^1 $.

If the $\theta^i$ exist, then it can be shown~\cite{DimMad96,MadMou96}
that the $\lambda_i$ must satisfy the equation
\be
2 \lambda_k \lambda_l P^{kl}{}_{ij} -
\lambda_k F^k{}_{ij} - K_{ij} = 0                                \label{eq:26}
\ee
with $F^k_{ij}$ and $K_{ij}$ complex numbers. Associated to this
equation there is a modified Yang-Baxter equation~\cite{MadMou96}.
The structure elements $C^i{}_{bc}$ are defined by the equation
\be
d\theta^i =
- {1\over 2} C^i{}_{jk} \theta^j \theta^k.                       \label{eq:27}
\ee
They are related to the coefficients of Equation~(\ref{eq:26}) by the
identity
\be
C^i{}_{jk} = F^i{}_{jk} - 2 \lambda_l P^{(li)}{}_{jk}.           \label{eq:28}
\ee
Consistent with Equation~(\ref{eq:25}), we shall impose the conditions
\be
P^{ij}{}_{kl} C^m{}_{ij} = C^m{}_{kl}, \qquad
P^{ij}{}_{kl} K_{ij} = K_{kl}.                                   \label{eq:29}
\ee

Using the Stehbein, we introduce now a connection and a torsion 2-form
\bea
D\theta^i & = & -\omega^i{}_{jk} \theta^j \otimes \theta^k,\\
\Theta^i & = & d \theta^i - \pi \circ D \theta^i,
\eea
as well as a metric
\be
g(\theta^i \otimes \theta^j ) = g^{ij}.
\ee
The coefficients $\omega^i{}_{jk}$ and $g^{ij}$ must lie in ${\cal A}$.
Since $g$ is ${\cal A}$-bilinear and because of the condition
(\ref{eq:20}) the  $g^{ij}$ must lie in the center ${\cal Z(A)}$. We
also write
\be
\sigma(\theta^i \otimes \theta^j ) =
S^{ij}{}_{kl} \theta^k \otimes \theta^l.
\ee
Then again by the ${\cal A}$-bilinearity of $\sigma$, the coefficients
$S^{ij}{}_{kl}$ lie also in ${\cal Z(A)}$. The condition (\ref{eq:12})
becomes
\be
(S^{ij}{}_{kl} + \delta^i_k \delta^j_l ) P^{kl}{}_{mn} = 0.     \label{eq:32}
\ee
Using this notation, the metric-compatibility of a connection $D$ is
expressed as
\be
\omega^i{}_{jk} + \omega_{kl}{}^m S^{il}{}_{jm} = 0.            \label{eq:33}
\ee
The condition that the connection be torsion-free is given by
\be
(\Omega^{i}{}_{jk} - \frac{1}{2}C^{i}{}_{jk})P^{jk}{}_{lm} = 0. \label{eq:34}
\ee
The curvature $\pi_{12}D^2 $ can be written in terms of the frame as
\be
\pi_{12}D^2 \theta^i =
- \frac{1}{2}R^{i}{}_{jkl} \theta^k \theta^l \otimes \theta^j  \label{eq:35-0}
\ee
and we have
\be
\mbox{Ric} (\theta^i) =
{1 \over 2} R^i{}_{jkl} \theta^k g(\theta^l \otimes \theta^j). \label{eq:35-1}
\ee
It is given by
\be
\mbox{Ric} \, (\theta^i) = R^i{}_j \theta^j.                   \label{eq:35-2}
\ee

\sect{The $h$-deformed quantum plane }

\subsect{ Linear connections }

The $h$-deformed quantum plane is an associative
algebra ${\cal A}$ generated by noncommuting elements (`coordinates')
$x$ and $y$ such that
\be
xy - yx = h y^2,                                                 \label{eq:35}
\ee
where $h$ is a deformation parameter.  The quantum group $GL_{h}(2)$ is
the symmetry group of the $h$-deformed plane as is $GL_{q}(2)$ for the
$q$-deformed quantum plane~\cite{DemManMukZhd90,Kup92}.
Let $T = \left( \begin{array}{cc}
                   A & B \\
                   C & D \end{array} \right) \in GL_{h}(2)$.
The commutation relations between the matrix elements of the quantum
group are given by
\bea
AB - BA & = & h\delta - h A^2, \nonumber \\
AC - CA & = &  h C^2,          \nonumber \\
AD - DA & = & h CD -h CA,      \nonumber \\
BC - CB & = & h CD + h AC,     \nonumber \\
BD - DB & = & h D^2 - h \delta,\nonumber \\
CD - DC & = & - h C^2,         \nonumber
\eea
where the quantum determinant
\be
\delta   = AD - CB - h CD = DA - CB - h CA
 \ee
is central.

$R$-matrix associated with this quantum group and which solves
the quantum Yang-Baxter equation
\be
\hat{R}_{12}\hat{R}_{23}\hat{R}_{12}
= \hat{R}_{23}\hat{R}_{12}\hat{R}_{23}                           \label{eq:38}
\ee
is given by
\be
\hat{R} = \left( \begin{array}{cccc}
             1  & -h & h & h^2 \\
             0  &  0  & 1 & h \\
             0  &  1 & 0 & -h  \\
             0  &  0 & 0 & 1
                \end{array}
              \right).                                          \label{eq:39}
\ee
The covariant differential calculus on the quantum plane can be
found~\cite{Agh93} by the method of Wess and Zumino~\cite{WesZum90}.
The results to be used in this work can be summarized as
follows.  For $x^i = (x, y)$ and $\xi^i = dx^i = (\xi, \eta )$ we have
\bea
&x^a x^b =        \hat{R}^{ab}{}_{cd}   x^c  x^d, \hspace{1cm}
&x^a \xi^b =      \hat{R}^{ab}{}_{cd} \xi^c  x^d,  \nonumber   \\
&\xi^a \xi^b =  - \hat{R}^{ab}{}_{cd} \xi^c\xi^d, \hspace{1cm}
&\partial_a x^b = \delta_a^{b}
                + \hat{R}^{bd}{}_{ac}x^c \partial_{d}.           \label{eq:40}
\eea
The second and third equations are written explicitly as
\be
\begin{array}{ll}
x \xi = \xi x - h \xi y +h \eta x + h^2 \eta y,
&x \eta = \eta x +h \eta y, \\
y \xi = \xi y - h \eta y,
&y \eta = \eta y,
\end{array}                                                     \label{eq:41}
\ee
and
\be
\xi^2 =  h \xi \eta,
\qquad \xi \eta = - \eta \xi,  \qquad
\eta^2 = 0.                                                     \label{eq:42}
\ee

Now, as in the $q$-deformed quantum plane~\cite{DubMadMasMou95}, the
action of a linear connection on the second equation of the above
relations generically results in the following relations
\be
\xi^a \otimes \xi^b = \hat{R}^{ab}{}_{cd}\sigma (\xi^c \otimes \xi^d ), \qquad
x^a D\xi^b = \hat{R}^{ab}{}_{cd} (D \xi^c ) x^d.                 \label{eq:43}
\ee
Since $\hat{R}^{-1} = \hat{R}$, the first equation is verified when
$\sigma $ transforms as $\hat{R} $, i.e.,
\bea
\sigma (\xi \otimes \xi ) & = &   \xi \otimes \xi - h \xi \otimes \eta
                            + h \eta \otimes \xi
                            + h^2 \eta \otimes \eta, \nonumber \\
\sigma (\xi \otimes \eta ) & = &   \eta  \otimes \xi
                            +  h \eta \otimes \eta, \\
\sigma (\eta \otimes \xi ) & = &   \xi \otimes \eta
                            - h \eta \otimes \eta, \nonumber \\
\sigma (\eta \otimes \eta ) & = &   \eta \otimes \eta.\nonumber  \label{eq:44}
\eea
The definition of the map $\sigma $ now is extended to the whole space
$\Omega^1 \otimes _{{\cal A}} \Omega^1 $ by the ${\cal A}$-linearity.
Then it is easy to see that $\sigma$ satisfies  Equation~(\ref{eq:12})
and
\be
\sigma ^2 = 1.
\ee
Also from the quantum Yang-Baxter equation of $\hat{R}$ in
Equation~(\ref{eq:38}) it follows that
\be
\sigma_{12}\sigma_{23}\sigma_{12}  = \sigma_{23}\sigma_{12}\sigma_{23}.
\ee

We are now in a position to exhibit the general linear connection
$D$ on $\Omega^1$, as in the $q$-deformed quantum plane. We introduce
the 1-form
\be
\kappa = x \eta - y\xi - h y \eta,
\ee
which is covariant under the action of $SL_{h}(2)$.
Then we have
\bea
&x \kappa = \kappa x,  \hspace{1cm} y\kappa = \kappa y, \nonumber \\
&\xi \kappa = - \kappa \xi,  \hspace{1cm} \eta \kappa = - \kappa \eta,
\eea
and
\be
\kappa^2 = 0 .                                                  \label{eq:49}
\ee
The $\kappa$ also obeys
\bea
\sigma (\xi \otimes \kappa ) & = & \kappa \otimes \xi, \hspace{1cm}
\sigma (\kappa \otimes \xi ) = \xi \otimes \kappa, \nonumber \\
\sigma (\eta \otimes \kappa )& = & \kappa \otimes \eta, \hspace{1cm}
\sigma (\kappa \otimes \eta ) = \eta \otimes \kappa, \\        \label{eq:50}
\sigma (\kappa \otimes \kappa ) & = & \kappa \otimes \kappa.\nonumber
\eea
A solution $D\xi^b $ to the second equation in (3.1.8) can be
immediately read off from the second equation of (3.1.5)
\be
D\xi^a = \rho (\xi^a \otimes \kappa + \kappa \otimes \xi^a ),
\label{eq:501}
\ee
where $\rho $ is a real  parameter.
Also, the first equation in (3.1.5) suggests  a solution of the form
\be
D\xi^a = \mu x^a \varpi,
\label{eq:502}
\ee
 where $\mu $ is a real parameter and $\varpi $ is any 1-form such that
$\pi \varpi = 0 $ and
\be
x^a \varpi = \varpi x^a.
\ee
 From Equation (3.1.13), it is easy to see that
$\varpi = \kappa \otimes \kappa $.
Then the general torsion-free linear connection $D$ is given by
\be
D\xi^a = \mu x^a \kappa \otimes \kappa +
 \rho (\xi^a \otimes \kappa + \kappa \otimes \xi^a ).
\label{eq:505}
\ee
This 2-parameter solution has been also found by Khorrami {\it et
al.}~\cite{KhoShaAgh97}. Now it is natural
to investigate other possible $\varpi$ in various cases. For this we
extend the $h$-deformed quantum plane to the two parameter
case~\cite{Agh93} on which the two parameter quantum group $GL_{h, h'}$
acts. In this case, we have the same equations as in (3.1.5) and (3.1.8)
with $\hat{R} $  replaced by
\be
\hat{R} = \left( \begin{array}{cccc}
             1  & -h' & h' & hh' \\
             0  &  0  & 1 & h \\
             0  &  1 & 0 & -h  \\
             0  &  0 & 0 & 1
                \end{array}
              \right).                              \label{eq:506}
\ee
A straightforward calculation yields 1-forms $\varpi $ as follows.
For $h' = nh $ $ (n = 2, 3, 4, \cdots ) $
\be
\varpi = y^{n-2} ( \kappa \otimes \eta + \eta \otimes \kappa ),
\ee
and for $h' = \frac{1}{2}nh $ $(n = 2, 3, 4, \cdots ) $
\be
\varpi = y^{n-2} \kappa \otimes \kappa
         - \frac{n-2}{2}hy^{n-1}\eta \otimes \kappa.
\ee
Then we have a 2-parameter family of torsion-free linear connections
 for $h' = nh $ $(n = 2, 3, 4, \cdots ) $
\be
D\xi^a =  \mu x^a ( y^{2n-2} \kappa \otimes \kappa
         - (n-1)hy^{2n-1} \eta \otimes \kappa )
      + \rho x^a y^{n-2}( \kappa \otimes \eta + \eta \otimes \kappa ),
\ee
and a 1-parameter family of  torsion-free linear connections
for  $h' = \frac{1}{2}nh $ $ (n = 3, 5, 7, \cdots ) $
\be
D\xi^a = \mu x^a ( y^{n-2} \kappa \otimes \kappa
         - \frac{n-2}{2}hy^{n-1}\eta \otimes \kappa ).
\ee
The supplementary condition $h' = nh $ or $h' = \frac{1}{2}nh $
corresponds to $p = q^n$ for the case~\cite{GeoMasWal96} of the
$q$-deformed quantum plane with the 2-parameter quantum group
$GL_{q, p}(2)$ where the linear connection is given by
$D\xi^a = \mu x^a x^{n-1}y^{n-1} \kappa \otimes \kappa$.

In the next Subsection, we shall concern the first term of the
2-parameter family of  torsion-free linear connections in
Equation~(\ref{eq:505})
\be
D\xi^a = \mu x^a \kappa \otimes \kappa                        \label{eq:508}
\ee
since it is compatible with a skew-symmetric nondegenerate bilinear map.
 From the linear connection  we have the curvature
\be
\pi _{12} D^2 \xi^a = - \Omega^a{}_{b}  \otimes \xi^b,
\ee
where the 2-form $\Omega^a{}_b $ is given by
\be
\Omega^a{}_b = 4 \mu \left(\begin{array}{cc}
                            xy & - x^2 + h xy \\
                           y^2 & -yx + hy^2
                           \end{array}\right) \xi\eta.
\ee
The 1-form $\kappa$ satisfies the equation $D^2 \kappa = 0$.

\subsect{ The symplectic 2-form }

In this Subsection, we shall use the expression `symplectic 2-form' and
`skew-symmetric metric' synonymously and denote it by $\Lambda $ since a
skew-symmetric metric on the $h$-deformed quantum plane resembles a
symplectic 2-form as in the ordinary geometry. A symmetric metric will
be denoted simply a metric as in ordinary geometry.

It has been shown that no metric can exist in the case of the $q$-deformed
quantum plane~\cite{DubMadMasMou95}.  However, the $h$-deformed quantum
plane has a better geometrical structure than the $q$-deformed quantum
plane and it does have a metric. In fact, the $h$-deformed quantum plane
has a symplectic 2-form, with which a metric can be associated.
The symplectic 2-form  of the $h$-deformed quantum plane is given
in a matrix form as
\be
\Lambda (\xi^a \otimes \xi^b ) \equiv
\Lambda^{ab} = \left( \begin{array}{cc}
                     h & 1 \\
                    -1 & 0 \end{array}
                \right).                                        \label{eq:54}
\ee
Now it is straightforward to show that in the particular case when the
covariant derivative is given by
$D\xi^a = \mu x^a \kappa \otimes \kappa$ we have
\be
( 1 \otimes \Lambda )D(\xi^a \otimes \xi^b ) =
 d \Lambda^{ab} = 0                                             \label{eq:55}
\ee
and for $\sigma_{23} = 1  \otimes \sigma$
\bea
(1 \otimes \Lambda )\sigma_{12}\sigma_{23}
              ( \xi \otimes \xi \otimes \xi ) &=& h \xi,  \nonumber \\
(1 \otimes \Lambda )\sigma_{12}\sigma_{23}
              ( \xi \otimes \xi \otimes \eta ) &=& h \eta,  \nonumber \\
(1 \otimes \Lambda )\sigma_{12}\sigma_{23}
              ( \xi \otimes \eta \otimes \xi ) &=&  \xi,  \nonumber \\
(1 \otimes \Lambda )\sigma_{12}\sigma_{23}
              ( \xi \otimes \eta \otimes \eta )& =&  \eta,   \\
(1 \otimes \Lambda )\sigma_{12}\sigma_{23}
              ( \eta \otimes \xi \otimes \xi ) &=& - \xi,  \nonumber \\
(1 \otimes \Lambda )\sigma_{12}\sigma_{23}
              ( \eta \otimes \xi \otimes \eta ) &=& - \eta,  \nonumber \\
(1 \otimes \Lambda )\sigma_{12}\sigma_{23}
              ( \eta \otimes \eta \otimes \xi ) &= & 0,  \nonumber \\
(1 \otimes \Lambda )\sigma_{12}\sigma_{23}
              ( \eta \otimes \eta \otimes \eta ) &=& 0.  \nonumber
\eea
 From these relations, it follows that the symplectic 2-form $\Lambda$
satisfies the compatibility condition in Equation~(\ref{eq:14}),
while the symplectic 2-form
$\Lambda = \left( \begin{array}{cc}
                     0 & 1 \\
                    -q & 0 \end{array}
                \right)
$ of the $q$-deformed quantum plane does not.

We define the $h$-deformed symplectic group by
\be
Sp_{h}(1) = \{T \in GL_{h}(2) \mid T \Lambda T^t = \Lambda \}.
\ee
Equivalently, the symplectic 2-form $\Lambda $ is
preserved under the action of $Sp_{h}(1)$, i.e.
\be
\Lambda ( \xi^{\prime a} \otimes \xi^{\prime b} )
 = \Lambda ( \xi^{a} \otimes \xi^{b} )
\ee
under the transformation
$\left( \begin{array}{c}
                     \xi^\prime \\
                     \eta^\prime  \end{array}
                \right) =
\left( \begin{array}{cc}
                     A & B \\
                     C & D \end{array}
                \right)
\left( \begin{array}{c}
                     \xi \\
                     \eta \end{array}
                \right)$.
 From  Equation~(\ref{eq:40}) it follows that $Sp_{h}(1) = SL_{h}(2)$,
which is consistent with the commutative limit when $h \rightarrow 0 $.

In ordinary symplectic geometry, metrics can be defined by a symplectic
2-form together with a complex structure and these all together result
in an Hermitian inner product. One can do the same in the $h$-deformed
quantum plane. We can define an ${\cal A}$-linear map
$J : \Omega^1 \rightarrow \Omega^1$ by
\be
J\xi = i \xi, \hspace{1cm} J \eta = -i \eta,
\ee
where $i = \sqrt{-1}$. The map $J$ satisfies $J^2 = -1 $ and can be
regarded as the complex structure of the $h$-deformed quantum plane.
Associated with the symplectic 2-form $\Lambda$, there is a metric $g$
satisfying the following relation, for $\xi, \eta \in \Omega^1 $,
\be
g(J \xi \otimes \eta ) = \Lambda ( \xi \otimes \eta ),
\ee
which can be written in matrix form as
\be
g = \left( \begin{array}{cc}
            -ih & -i \\
            -i & 0
           \end{array}
   \right).
\ee
On the other hand, there is another metric $g^\prime$  defined by
\be
g^\prime(\xi \otimes \eta) = \Lambda (\xi \otimes J \eta),
\ee
which is, in matrix form,
\be
g^\prime = \left( \begin{array}{cc}
            ih & -i \\
            -i & 0
           \end{array}
   \right).
\ee
The two metrics are related by the condition
\be
g(J\xi \otimes J\eta ) = g^\prime(\xi \otimes \eta )
\ee
for any $\xi, \eta \in \Omega^1$ and agree when $h \rightarrow 0$.
These metrics, however, are not compatible with the linear connection
$D$.  In fact there is no metric on the $h$-deformed plane with respect
to which $D$ is compatible. Such a metric can be found however if we
extend the $h$-deformed quantum plane as in the next section.  In order
to compare them with the commutative-limit case let us define
\be
\vartheta^1 = \frac{1}{\sqrt 2 }(\xi + i \eta ),
\hspace{1cm} \vartheta^2 = \frac{i}{\sqrt 2 }(\xi - i \eta ).
\ee
Then it is easy to see that
\be
J\vartheta^1 = \vartheta^2,
\hspace{1cm}
J\vartheta^2 = - \vartheta^1.
\ee
With respect to $\{ \vartheta^1, \vartheta^2 \}$, the two metrics
$g,\,\, g^\prime$, and the symplectic 2-form $\Lambda $
can be expressed as follows
\be
g = \left( \begin{array}{cc}
            1-\frac{ih}{2} & \frac{h}{2} \\
            \frac{h}{2} & 1 + \frac{ih}{2}
           \end{array}
   \right),
\hspace{0.3cm}
g^\prime = \left( \begin{array}{cc}
            1 + \frac{ih}{2} & - \frac{h}{2} \\
            -\frac{h}{2} & 1 - \frac{ih}{2}
           \end{array}
   \right),
\hspace{0.3cm}
\Lambda = \left( \begin{array}{cc}
            \frac{h}{2} & 1 + \frac{ih}{2} \\
           -1 + \frac{ih}{2} & - \frac{h}{2}
           \end{array}
   \right).
\ee

If we define
\be
 H = g^\prime + i\Lambda,
\ee
the map $H$ goes over to the usual Hermitian inner product on
the complex 2-plane ${\mathbb C}^2 $ in the commutative limit.
Moreover, it is interesting to see that
if we let $\eta_{a} = (-\eta, \xi + h \eta )$ and
define the skew derivative $\partial_a$ by
\be
\partial_a f = \Lambda (\eta_a \otimes df ),
\ee
then the $\partial_a$ satisfy the second equation of~(\ref{eq:40}) given
by Wess and Zumino~\cite{WesZum90}. Thus the skew derivatives
$\partial_a$ arise as Hamiltonian vector fields in the $h$-deformed
quantum plane.  This is not the case for the $q$-deformed quantum plane.
In fact, if $x^a \xi^b = \hat{R}^{ab}{}_{cd}\xi^c x^d $, there should be
elements $\eta_a \in \Omega^1$ such that
\be
\eta_a x^b = \hat{R}^{bc}{}_{ad} x^d \eta_{c}
\ee
for the skew derivatives to be induced from the symplectic 2-form
$\Lambda $ as above.  However, there are no such $\eta_a$ in the case of
the $q$-deformed quantum plane.

\sect{The extended $h$-deformed plane}

\subsect{Linear connections}

The extended $h$-deformed quantum plane is an associative algebra
${\cal A}$ generated by $x, y, x^{-1}, y^{-1}$ satisfying
Equation~(\ref{eq:35}). The extended $h$-deformed plane is also more
interesting than the extended $q$-deformed one from point of view of
geometry since the metric and linear connection it supports
have an interesting commutative limit.

If ${\cal A}$ is a unital $*$-algebra and $x$ and $y$ are Hermitian elements,
then $h \in i {\mathbb R}$.  Equation~(\ref{eq:35}) can be
written as $[x, y^{-1}] = - h$.  In this form we see that the algebra
has the structure of the Heisenberg algebra with the parameter $h$
playing the role of $\hbar$ but the differential calculus (\ref{eq:41})
is not `natural' from this point of view.  If we introduce
\be
u = x y^{-1} + {1\over 2} h, \qquad v = y^{-2}                     \label{71}
\ee
then the commutation relation becomes
\be
[u,v] = - 2 h v.                                                \label{eq:72}
\ee
This choice of generators is useful in studying the commutative limit.
If $x$ and $y$ are Hermitian, then so are $u$ and $v$.

We can write (\ref{eq:72}) also as $[u, (1/2)\log v] = - h$ if we
introduce  the formal element $\log v$. We see then that
$$
x^\prime = u, \qquad y^\prime = {2 \over \log v}
$$
satisfy also the commutation relations (\ref{eq:35}). The algebra cannot then
be uniquely defined by the commutation relations.  In fact von~Neumann
proved that only by using additional topological conditions could one
deduce the uniqueness of the representation of the Heisenberg
commutation relations.

A (real) frame can be written in terms of the generators in
Equation~(\ref{eq:35}) as
\be
\theta^1 = y \xi - (x - h y) \eta, \qquad \theta^2 = 2 y^{-1} \eta
\ee
and in terms of the generators in Equation~ (\ref{eq:72}) as
\be
\theta^1 = v^{-1} du, \qquad \theta^2 = - v^{-1} dv.            \label{eq:74}
\ee
The original basis $(\xi, \eta)$ can be written in terms of the $\theta^a$ as
\be
2 \xi = 2 y^{-1} \theta^1 + x \theta^2, \qquad 2 \eta = y  \theta^2.
\ee
The $\theta^a$ satisfy the commutation relations~(\ref{eq:20}) as well
as the relations
\be
(\theta^1)^2 = 0, \qquad (\theta^2)^2 = 0, \qquad
\theta^1 \theta^2 + \theta^2 \theta^1 = 0.                       \label{eq:77}
\ee
 From Equation~(\ref{eq:23}) we see then that $P^{ab}{}_{cd}$ is given by
\be
P^{ab}{}_{cd} = {1\over 2} (\delta^a_c \delta^b_d - \delta^b_c \delta^a_d)
                                                                 \label{eq:78}
\ee
and therefore, from Equation~(\ref{eq:28}), we have
$C^a{}_{bc} = F^a{}_{bc}$. In particular the $C^a{}_{bc}$ are real
numbers.  The differentials $d\theta^a$ are given by
Equation~(\ref{eq:27}) with
\be
C^1{}_{12} = - C^1{}_{21} = 1, \qquad C^2{}_{ab} = 0.
\ee

If we introduce the derivations $e_a = \ad \lambda_a$ with
\be
\lambda_1 = {1 \over 2h} y^{-2} = {1 \over 2h} v, \qquad
\lambda_2 = {1 \over 2h} x y^{-1} + {1\over 4} = {1 \over 2h} u. \label{eq:81}
\ee
we see that Equation~(\ref{eq:19}) is satisfied. We can conclude from
Equation~(\ref{eq:26}) that the $\lambda_a$ must form a (real) Lie algebra.
We have then from Equation~(\ref{eq:72})
\be
[\lambda_1, \lambda_2] = \lambda_1.                              \label{eq:82}
\ee
The $\lambda_a$ form a solvable Lie algebra.

The `Dirac operator'~\cite{Con94} in Equation~(\ref{eq:21}) is given by
\be
\theta = -{1\over 2h} y^{-1} \xi - {1\over 2h} (x -hy) y^{-2} \eta
= -{1\over 2h} (du - u v^{-1} dv).                               \label{eq:80}
\ee
A straightforward calculation yields $d\theta + \theta^2 = 0$.

We introduce a metric and we set $g(\theta^a \otimes \theta^b) = g^{ab}$.
 From the bilinearity condition on $g$ and the relations (\ref{eq:35}) we
see that the $g^{ab}$ must be complex numbers. If we wish the metric to
be real then the $g^{ab}$ must be real numbers. By a trivial change of
basis we can suppose that $g^{ab} = \delta^{ab}$.  We have then in terms
of the generators $x$ and $y$
\be
\begin{array}{ll}
g(\xi \otimes \xi) = y^{-2} +  x^2/4, &g(\xi \otimes \eta) =  xy/4, \\
g(\eta \otimes \xi) =  yx/4,          &g(\eta \otimes \eta) =  y^2/4
\end{array}                                                      \label{eq:83}
\ee
and in terms of the generators $u$ and $v$
\be
\begin{array}{ll}
g(du \otimes du) = v^2, &g(du \otimes dv) = 0, \\
g(dv \otimes du) = 0,   &g(dv \otimes dv) = v^2.
\end{array}
\ee
A flat metric-compatible linear connection is given by
\be
D\theta^a = 0.
\ee
It has torsion, given by
\be
\Theta^1 = - \theta^1 \theta^2, \qquad \Theta^2 = 0.
\ee
The unique torsion-free, metric-compatible linear connection is given by
\be
D\theta^1 = - \theta^1 \otimes \theta^2, \qquad
D\theta^2 = \theta^1 \otimes \theta^1.
\ee
This $D$ is also compatible with the symplectic 2-form $\Lambda$ given
in Equation~(\ref{eq:54}): $D\Lambda = 0$.

The curvature map defined by Equation~(\ref{eq:35-0}) becomes
\be
\pi_{12}D^2 \theta^1 =  \theta^1 \theta^2 \otimes \theta^2, \qquad
\pi_{12}D^2 \theta^2 =  - \theta^1 \theta^2 \otimes \theta^1.
\ee
If one sets as usual $R_{abcd} = g_{ae} R^e{}_{bcd}$ then one finds that
the Gaussian curvature is given by
\be
R_{1212} = - 1.                                                 \label{eq:89}
\ee
The coefficients $R_{abcd}$ satisfy the usual symmetries of the
coefficients of a Riemann curvature tensor. The coefficients of the
Ricci map are given by
\be
R^a{}_b = \delta^a_b.
\ee

We choose now $n = 3$. Then it is of interest to introduce a third
(Hermitian) element
\be
w = - {1\over 2}(u^2 - 2hu + 1 + 2h^2 )v^{-1}                   \label{eq:91}
\ee
of the algebra ${\cal A}$ and define
\be
\lambda_3 = {1\over 2h} w.
\ee
The $\lambda_i = (\lambda_a, \lambda_3)$ still form a Lie algebra
\be
[\lambda_1, \lambda_2] = \lambda_1, \qquad
[\lambda_2, \lambda_3] = \lambda_3, \qquad
[\lambda_3, \lambda_1] = \lambda_2.                             \label{eq:93}
\ee
A straightforward calculation yields
\be
\begin{array}{lll}
e_1 u = v,   &e_1 v = 0,   &e_1 w = - u, \\
e_2 u = 0,   &e_2 v = - v, &e_2 w = w,   \\
e_3 u = - w, &e_3 v = u,   &e_3 w = 0.
\end{array}                                                      \label{eq:94}
\ee
The $e_i$ satisfy the same commutation relations
\be
[e_1, e_2] = e_1, \qquad [e_2, e_3] = e_3, \qquad [e_3, e_1] = e_2.
                                                                 \label{eq:95}
\ee
as the $\lambda_i$.  They are real in the sense that the derivation
$e_i f$ of an Hermitian element $f$ is Hermitian.

The Lie algebra in Equation~(\ref{eq:95}) is a real form of
$SL(2, {\mathbb C})$, different from the Lie algebra of $SO_3$.  We have
found a frame with 2 generators since the Poincar\'e half-plane is a
parallelizable manifold and the module of 1-forms is a free (left or
right) module. This is not so in the case of the 2-sphere~\cite{Mad95};
the module of 1-forms in this case is a nontrivial submodule of a free
module of rank 3. The Lie algebra of Killing vector fields of the
Poincar\'e half-plane and the sphere are different real realizations of
$SL(2, {\mathbb C})$.

A differential calculus can be defined using the three 1-forms
$\theta^i$ dual to the derivations $e_i$. An analogous situation was
discussed in the case of the $q$-deformed plane~\cite{DimMad96}.
 From Equation~(\ref{eq:94}) we conclude that
\be
du = v \theta^1 - w \theta^3, \qquad
dv = - v \theta^2 + u \theta^3,                                  \label{eq:97}
\ee
to which we can add
\be
d w = - u \theta^1 + w \theta^2.                                 \label{eq:98}
\ee
The second of the Equations~(\ref{eq:97}) is a trivial consequence of the
commutation relations, obtained by equating the differential of both
sides of Equation~(\ref{eq:72}).  The previous differential calculus
with two generators is obtained formally by setting $\theta^3 = 0$ in
Equation~(\ref{eq:97}). The commutation relations which define the
module structure of $\Omega^1$ are obtained from Equation~(\ref{eq:97}):
\be
\begin{array}{ll}
udu -duu = - 2hdu - 4hw\theta^3, &vdu -duv = 2hu \theta^3,   \\
udv -dvu = - 2hdv + 2hu\theta^3, &vdv -dvv = 2hv \theta^3.
\end{array}
\ee
Apart from the trivial relation which follows from the commutation
relations these equations contain two cubic relations
\be
\begin{array}{l}
uvdv - udvv = vduv - duv^2, \\
vudu - vduu + 2hvdu + 2vdvw - 2dvvw = 0.
\end{array}
\ee
Provided that $h \neq 0$ the system Equation~(\ref{eq:97}), (\ref{eq:98})
can be inverted to obtain equations for the $\theta^i$ in terms of $du$,
$dv$ and $dw$:
\be
\theta^1 = {1\over 2h} u^{-1} [w, du], \qquad
\theta^2 = {1\over 2h} v^{-1} [u, dv], \qquad
\theta^3 = {1\over 2h} u^{-1} [v, du].                         \label{eq:101}
\ee
This differential calculus has fewer relations than the one defined
above. It lies between the one defined by the relations (\ref{eq:77}) and the
universal differential calculus, which has a free algebra of forms with
no relations.

One can define a Lie derivative in noncommutative geometry exactly as
the ordinary Lie derivative is defined in ordinary geometry. If $\xi$
is a form and $i_X$ is the interior product then the Lie derivative
$L_X \xi$ of $\xi$ is given by
\be
L_X \xi = d i_X \xi + i_X d \xi
\ee
A Killing derivation~\cite{Mad95} is a derivation $X$ such that the
Lie derivative $L_X$ of the metric $g$ vanishes: $L_X g = 0$. Denote
$L_a$ the Lie derivative with respect to the derivation $e_a$. Then it
is easy to see that
\be
\begin{array}{ll}
L_1 \theta^1 = - \theta^2,          &L_1 \theta^2 = 0,\\
L_2 \theta^1 = + \theta^1,          &L_2 \theta^2 = 0,\\
L_3 \theta^1 = - v^{-1} w \theta^2, &
L_3 \theta^2 = - v^{-1} u \theta^2 - \theta^1.
\end{array}
\ee
 From these formulae one can calculate the Lie derivative of the metric:
\be
\begin{array}{l}
L_1 g = - (\theta^1 \otimes \theta^2 + \theta^2 \otimes \theta^1),\\
L_2 g = \phantom{+} 2 \theta^1 \otimes \theta^1,\\
L_3 g = - (1 + v^{-1} w) (\theta^1 \otimes \theta^2 +
\theta^2 \otimes \theta^1) - 2 v^{-1} u \theta^2 \otimes \theta^2.
\end{array}                                                   \label{eq:101a}
\ee
That is, none of the derivations $e_i$ is a Killing derivation.

\subsect{ The commutative limit of the extended plane }

It is interesting to study the structure of the extended $h$-deformed
quantum plane in the commutative limit.  In terms of the commutative
limits $\tilde u$, $\tilde v$ of the generators $u$, $v$ of the algebra
${\cal A}$ and the corresponding commutative limit $\tilde \theta^a$ of
the frame, the metric is given by the line element
\be
ds^2 = (\tilde \theta^1)^2 + (\tilde \theta^2)^2 =
\tilde v^{-2} (d \tilde u^2 + d \tilde v^2).                  \label{eq:102}
\ee
This is the metric of the Poincar\'e half-plane. The algebra
${\cal A}$ with the differential calculus defined by the relations
(\ref{eq:77}) can be considered then as a noncommutative deformation of the
Poincar\'e half-plane.

The derivations $e_i$ define, in the commutative limit, 3 vector fields
\be
X_i = \lim_{h \rightarrow 0} e_i.                               \label{eq:103}
\ee
If we define $\tilde w$ to be the commutative limit of $w$ then
\be
X_1 = \tilde v \partial_{\tilde u}, \qquad
X_2 = - \tilde v \partial_{\tilde v}, \qquad
X_3 = - \tilde w \partial_{\tilde u} + \tilde u \partial_{\tilde v}.
                                                                \label{eq:104}
\ee
By construction these vector fields form a Lie algebra with the same
commutation relations as the $e_i$.  By the Equation~(\ref{eq:101a}) of
the previous section we see that the $X_i$ cannot be Killing vector
fields. There is in fact no reason for them to be so.  The Poincar\'e
half-plane has however three Killing vector fields $X^\prime_i$, given
by
\be
X^\prime_1 = \partial_{\tilde u}, \qquad
X^\prime_2 = \tilde u \partial_{\tilde u} +
\tilde v \partial_{\tilde v}, \qquad
X^\prime_3 = {1\over 2} (\tilde v^2 - \tilde u^2) \partial_{\tilde u}
-  \tilde u \tilde v \partial_{\tilde v}.                       \label{eq:106}
\ee

Define a map $\phi$ of the Poincar\'e half-plane into itself by
\be
\phi(\tilde u) = \tilde u^\prime = \tilde u \tilde v^{-1}, \qquad
\phi(\tilde v) = \tilde v^\prime = \tilde v^{-1}.
\ee
This is a regular diffeomorphism. Indeed
\be
\phi^2 = \phi \circ \phi = 1.
\ee
In the spirit of noncommutative geometry we consider $\tilde u$ and
$\tilde v$ as generators of the algebra of functions on the Poincar\'e
half-plane. In ordinary differential geometry a map $\phi$ of the
manifold induces a map $\phi^*$ of the algebra of differential
forms and a map Let $\phi_*$ of the vector fields. Since we shall
not have occasion to refer to the manifold as such we use the notation
$\phi$ to designate the restriction of $\phi^*$ to the algebra of
functions. Since we have
\be
\phi_*\partial_{\tilde u} = \partial_{\tilde u^\prime} =
\tilde v \partial_{\tilde u}, \qquad
\phi_*\partial_{\tilde v} = \partial_{\tilde v^\prime} =
- \tilde u \tilde v \partial_{\tilde u}
- \tilde v^2 \partial_{\tilde v}                                \label{eq:107}
\ee
it is easy to see that
\be
\phi_* X_i = X^\prime_i.
\ee
The commutative limit of the derivations which defined the differential
calculus are related to the Killing vector fields then in a simple way.
We have not succeeded in constructing derivations of the algebra whose
commutative limits are the Killing vector fields $X^\prime_i$.  The
limit $h \rightarrow 0$ is a rather singular limit and it need not be
true that an arbitrary vector field on the Poincar\'e half-plane is the
limit of a derivation. The action of $\phi^*$ on the frame is given by
\be
\phi^* \tilde \theta^1 = \tilde v^{\prime -1} d \tilde u^\prime
= \tilde v \tilde \theta^1 + \tilde u \tilde \theta^2, \qquad
\phi^* \tilde \theta^2 = - \tilde v^{\prime -1} d \tilde v^\prime
= - \tilde \theta^2.
\ee
The vector fields $X_i$ are Killing with respect to the metric
\be
ds^2 = (\phi^* \tilde \theta^1)^2 +
       (\phi^* \tilde \theta^2)^2.                    \label{eq:108}
\ee

The map $\phi$ can also be considered as a change of coordinates.
In this case Equations~(\ref{eq:102}) and (\ref{eq:108}) describe
the same line element in different coordinates systems. The components
$X^a_i$ of the vector fields $X_i$ coincide with the components of the
Killing vector fields in the new coordinate system and the components
$X^{\prime a}_i$ of the vector fields $X^\prime_i$ coincide with the
components of the Killing vector fields in the old coordinate system.
We have not really understood the role of the map $\phi$ nor why
it appears. We constructed the algebra ${\cal A}$ using generators and
relations.  This is the noncommutative version of the method of defining
a curved manifold as an embedding in a higher-dimensional flat euclidean
space.  It is known~\cite{Gam91, MleSte93} that the Poincar\'e
half-plane cannot even be immersed in ${\mathbb R}^3$.  This fact might
somehow also be connected with the existence of the map $\phi$.

The commutation relations (\ref{eq:72}) define on the Poincar\'e
half-plane a Poisson structure
\be
\{\tilde u, \tilde v\} = - 2 \tilde v.
\ee
Since the map $\phi$ is not a symplectomorphism it cannot be `lifted' to
a morphism of the algebra ${\cal A}$.  There should be a
relation~\cite{Mad97} between the Poisson structure and the Riemann
curvature. It is not evident from the present example however what this
relation could be. The Poincar\'e half-plane has been used as an example
of a classical and quantum phase space and as such has many interesting
properties.  For a discussion of this and reference to the previous
literature we refer to Emch {\it et al.}~\cite{EmcNarThiSew94}. The
relation between the algebra of a free quantum particle on the
Poincar\'e half-plane and the $h$-deformed algebra we have used has yet
to be investigated.

\sect{ Conclusion }

The $h$-deformed quantum plane seems to have more geometrical structures
than the $q$-deformed one.  In the $h$-deformed quantum plane, there is
a 2-parameter family of torsion-free linear connections.
The existence of a 2-parameter family of
torsion-free linear connections is  quite general even within
the set of 2-parameter $h$-deformed quantum planes.
Moreover, there is a skew-symmetric non-degenerate bilinear map with
which a 1-parameter sub-family of linear connections are compatible. 
The skew-symmetric  map plays an important role.  It
resembles the symplectic 2-form and makes the linear connections
symplectic.  Moreover, it is interesting that the skew-symmetric map
induces skew derivatives in the $h$-deformed quantum plane.  We can also
define a complex structure on the $h$-deformed quantum plane and
construct a metric using this structure together with the skew-symmetric
map as in the ordinary symplectic geometry~\cite{AebBorKalLeuRei92}.
However, it should be stressed that the metric is not compatible with
the  linear connections.  A similar
construction is not possible in the case of the $q$-deformed quantum
plane.

The geometry of the Poincar\'e half-plane can be completely globally
defined by the action of the $SL(2,{\mathbb R})$ group whose Lie algebra
is given by the Killing vectors. From this point of view a complete
classification of all Poisson structures on the Poincar\'e (Lobachevsky)
half-plane as well as their possible `quantum' deformations has been
given by Leitenberger~\cite{FadResTak89, Lei96}. We have analysed in
detail the extended $h$-deformed quantum plane as a noncommutative
version of the Poincar\'e half-plane; the roles of the derivations and
the Stehbein are explicitly investigated.

\section*{Acknowledgments}

This work was supported by Ministry of Education, Project No.\
BSRI-96-2414 and a grant from TGRC-KOSEF 1996. One of the authors (SC)
is also financially supported by the 1996 Research Program of Semyung
University and would like to thank all colleagues of LPTHE at
Universit\'{e} de Paris-Sud for their hospitality and helpful
discussions during his visit.  Another author (JM) would like to thank
the Ludwig-Maximilian Universit\"at for its hospitality and financial
support.  He is also grateful to A. Chakrabarty and T. Friedrich for
useful discussions. Finally the authors would like to thank M.  Khorrami
for informing us of his work~\cite{KhoShaAgh97} on the geometry of
$h$-deformations.


\begin{thebibliography}{99}
\parskip 5pt plus2pt minus2pt
\tolerance=1000

\bibitem{AbdChaCha97}
B. Abdesslam, A. Chakrabarti, R. Chakrabarti,
``General construction of nonstandard $R_{h}$-matrices as
contraction limits of $R_{q}$-matrices'',  q-alg/9706033.

\bibitem{AghKhoSha95}
A. Aghamohammadi, M. Khorrami, A. Shariati,
``$h$-deformation as a contraction of $q$-deformation'',
J. Phys. A {\bf  28}, (1995) L225.

\bibitem{AebBorKalLeuRei92}
B. Aebischer, M. Borer M. M. K\"alin, Ch. Leuenberger, H.M. Reimann,
``Symplectic geometry'', Birkh\"auser Verlag, Basel (1992).

\bibitem{Agh93}
A. Aghamohammadi,
``The two-parametric  extension of  $h$ deformation of $GL(2)$
and the differential calculus on its quantum plane'',
Mod. Phys. Lett. A {\bf 8} (1993) 2607.

\bibitem{ChaPre94}
V. Chari, A. Pressley, ``A Guide to Quantum Groups'',
Cambridge University Press, 1994,

\bibitem{Con94}
A. Connes, ``Noncommutative Geometry'',
Academic Press, 1994.

\bibitem{DemManMukZhd90}
E.E. Demidov, Yu.I. Manin, E.E. Mukhin, D.V.Zhdanovich
``Non-standard quantum deformations of $GL(n)$
and constant solutions of the Yang-Baxter equation'',
Prog. Theor. Phys. (Suppl.)  No. 102, (1990) 203.

\bibitem{DimMad96}
A. Dimakis, J. Madore, ``Differential Calculi and Linear Connections'',
J. Math. Phys. {\bf 37} (1996) 4647.

\bibitem{DubMadMasMou95}
M. Dubois-Violette, J. Madore, T. Masson, J. Mourad,
``Linear Connections on the Quantum Plane'',
Lett. Math. Phys. {\bf 35} (1995) 351.

\bibitem{DubMadMasMou96}
M. Dubois-Violette, J. Madore, T. Masson, J. Mourad,
``On Curvature in Noncommutative Geometry'',
J. Math. Phys. {\bf 37} (1996) 4089.

\bibitem{EmcNarThiSew94}
G.G. Emch, H. Narnhofer, W. Thirring and G.L. Sewell,
``Anosov actions on non-commutative algebras''
J. Math. Phys. {\bf 35} (1994) 5582.

\bibitem{FadResTak89}
L. D. Faddeev, N. Y. Reshetikhin and L. A. Takhtajan,
``Quantization of Lie Groups and Lie Algebras'',
Algebra i Analysis, {\bf 1} (1989), 178;
translation: Leningrad Math. J. {\bf 1} (1990) 193.

\bibitem{Gam91}
R.V. Gamkrelidze (Ed.) ``Geometry I'' , Encyclop\ae dia of Mathematical
Sciences {\bf 28}, Springer Verlag (1991).

\bibitem{GeoMasWal96}
Y. Georgelin, T. Masson, J-C Wallet,
``Linear connections on the two parameter quantum plane'',
Rev. Math. Phys. {\bf 8} (1996) 1055.

\bibitem{KhoShaAgh97}
M. Khorrami, A. Shariati, A. Aghamohammadi,
``$SL_{h}(2) $-Symmetric Torsionless Connections'',
Lett. Math. Phys. {\bf 40} (1997) 95

\bibitem{Kup92}
B.A. Kupershmit, ``The quantum group $GL_{h}(2)$'',
J. Phys. A {\bf 25} (1992) L1239.

\bibitem{Lei96}
F. Leitenberger, ``Quantum Lobachevsky Planes'',
J. Math. Phys. {\bf 37} (1996) 3131.

\bibitem{Mad95}
J. Madore, ``An Introduction to Noncommutative Differential Geometry
and its Physical Applications'', Cambridge University Press, 1995.

\bibitem{Mad97}
J. Madore, ``On Poisson Structure and Curvature'',
Preprint, LPTHE Orsay 97/25, gr-qc/9705083

\bibitem{MadMou96}
J. Madore, J. Mourad, ``Quantum Space-Time and Classical Gravity'',
J. Math. Phys. (to appear), gr-qc/9607060.

\bibitem{Maj95}
S. Majid, ``Foundations of Quantum Group Theory'',
Cambridge University Press, 1995.

\bibitem{Man88}
Yu.I. Manin, ``Quantum groups and Noncommutative geometry'',
Centre de Recherches Math\'ematiques, Montr\'eal (1988).

\bibitem{Man91}
Yu.I. Manin, ``Topics in Noncommutative Geometry'',
Princeton Uiversity Press, Princeton (1991).

\bibitem{MleSte93}
M. Mleko, I. Sterling, ``Application of soliton theory to the
construction of pseudospherical surfaces in ${\mathbb R}^3$'',
Annals of Global Analysis and Geometry, {\bf 11} (1993) 65.

\bibitem{Mou95}
J. Mourad, ``Linear connections in non-commutative geometry'',
Class. Quant. Grav. {\bf 12} (1995) 965.

\bibitem{PusWor89}
W. Pusz and S.L. Woronowicz, ``Twisted Second Quantization'',
Rep. on Math. Phys. {\bf 27} (1989) 231.

\bibitem{WesZum90}
J. Wess, B. Zumino,
``Covariant differential calculus on the quantum hyperplane'',
Nucl. Phys. B (Proc. Suppl.) {\bf 18} (1990) 302.

\end{thebibliography}
\end{document}